\documentclass{amsart}
\usepackage{amsrefs}
\usepackage[pdftex,colorlinks]{hyperref}

\newtheorem*{stratification-conjecture}{Rigidity Conjecture}

\newcommand\defn[1]{\textbf{#1}}

\begin{document}

\title{The rigidity conjecture}
\author{Marco Martens}
\author{Liviana Palmisano}
\author{Bj\"orn Winckler}
\address{Institute for Mathematical Sciences, Stony Brook University, Stony
    Brook NY 11794-3660, USA}
\address{School of Mathematics, University of Bristol, Bristol BS8 1TW, UK}
\email{marco@math.stonybrook.edu}
\email{liviana.palmisano@gmail.com}
\email{bjorn.winckler@gmail.com}
\date{\today}

\begin{abstract}
A central question in dynamics is whether the topology of a system determines
its geometry.
This is known as rigidity.
Under mild topological conditions rigidity holds for many classical cases,
including: Kleinian groups, circle diffeomorphisms, unimodal interval maps,
critical circle maps, and circle maps with a break point.
More recent developments show that under similar topological conditions,
rigidity does not hold for slightly more general systems.
In this paper we state a conjecture which describes how topological classes are
organized into rigidity classes.
\end{abstract}

\maketitle

\section{Introduction}

One of the aims of dynamics is to understand whether two dynamical systems
are ``topologically'' the same.
This is determined by the existence of a homeomorphism which conjugates the two
systems.
A related question is then to ask when two systems are ``geometrically'' the same.
That is, when is the conjugacy differentiable?

This geometrical equivalence question has been studied in the last forty years
in the case of circle diffeomorphisms, unimodal maps, critical circle maps,
etc.\ (see Example~\ref{standard}).
It turns out that, under mild topological restrictions, the conjugacy between
two systems is differentiable as soon as it exists.
In other words, the topology of a system determines its geometry.
This is called the rigidity phenomenon.

One cannot expect the rigidity phenomenon in all generality. The mild topological restrictions are essential. For example, there is no rigidity in the context of circle diffeomorphisms when the rotation number is of strongly unbounded type, \cites{A61, H79}. We will discuss the rigidity phenomenon only in the situation of bounded combinatorics.
This is done with the purpose of stressing the fact that even in
this simplest situation the rigidity phenomenon is more intricate than the classical
case where ``topology determines geometry.''

Only in the last few years, further studies about the geometry of dynamical
systems with bounded combinatorics have revealed classes for which the rigidity
phenomenon does not hold.
Non-rigidity occurs in natural classes of dynamical systems, such as: circle
maps with a flat interval, Lorenz maps in one dimension and in H\'enon maps in
two dimensions.
The geometrical equivalence of these systems is not solely determined by their
topology.
However, the rigidity phenomenon does not break down completely.
Instead, the geometrical equivalence classes, called rigidity classes, are well
organized inside the topological ones.
The observed structures are
\begin{itemize}
 \item foliations by rigidity classes,
 \item the coexistence phenomenon,
 \item probabilistic rigidity.
\end{itemize}
These notions are described in more detail in Sections \ref{rigconj}
and~\ref{examples}.

The above examples and the structures that they revealed are what urged us to
come up with a conjecture which describes the relation between the topological
and geometric properties of a system.
In Section~\ref{rigconj} we discuss the resulting Rigidity Conjecture and in
Section~\ref{examples} we give examples supporting it.

\section{The Rigidity Conjecture}\label{rigconj}

In this section we present the basic notions needed to state the Rigidity
Conjecture.
The aim is to determine the geometry of the attractor of a system.
The systems are smooth maps on manifold and the attractors are attractors in
the sense of Milnor \cite{Milnor}.

Two maps are in same \defn{topological class} if they are conjugated on their
attractors.
Similarly, two maps are in same \defn{rigidity class} if they are
$C^{1+\alpha}$--conjugated on their attractors, for some $\alpha>0$.
A third notion of equivalence is given by so-called probabilistic rigidity.
An attractor carries a dynamically relevant measure and we say that
two maps are in the same \defn{probabilistic rigidity class} if the
conjugacy is $C^{1+\alpha}$ almost everywhere with respect to this measure,
for some $\alpha > 0$.
The topological class determines the topological properties of the attractor,
whereas the rigidity class determines the attractor's geometrical properties.

We restrict our discussion to topological classes which are of
\defn{bounded combinatorics}.
This topological property is well understood for one-dimensional systems and for
infinitely renormalizable H\'enon maps.
For example, in case of circle diffeomorphisms bounded combinatorics is the same as
saying that the rotation number is of bounded type.
However, the topology of two and higher dimensional systems is still in the very beginning of its development. Part of the study of the rigidity phenomenon is to describe the topological restrictions needed for rigidity. At this moment our understanding of the topology of higher dimensional systems is too rudimentary to anticipate the general condition needed for rigidity. These topological restrictions will have the nature of being bounded.

Finally, a \defn{stratification} of a topological class is a partition of the
topological class into finite codimension submanifolds.
The submanifolds can have different codimensions.
Some of them can form a foliation.

\begin{stratification-conjecture}
The topological class is a finite codimension manifold which is stratified by
probabilistic rigidity classes.
\end{stratification-conjecture}

The heuristic reasoning behind the conjecture is as follows.
The rigidity classes determine the geometrical properties of the attractor on
a small scale.
The tool for studying small-scale properties of an attractor is
renormalization.
Renormalization allows us to zoom in on any point of the
attractor.\footnote{Usually the renormalization is done around one specific
point.  Nevertheless the same procedure can be used for any point of the
attractor.}
If the renormalizations of two maps at any point of the attractor converge
exponentially fast then the conjugacy is $C^{1+\alpha}$, with $\alpha>0$.
For the one-dimensional setting see \cites{dFdM, KT} and for the
conservative H\'enon case see \cite{GJM}.
Two maps are in the same rigidity class if and only if the renormalizations at
any point of the attractor converge exponentially fast. This relies on the restriction of bounded combinatorics \cites{dFdM, dFdM1, Av13, KK14, KT1, KT}.
Roughly speaking, the rigidity classes have finite codimension because the
derivatives of the renormalization are compact operators.
The rigidity classes are then determined by the asymptotic behavior of the
renormalizations, which in turn strongly depend on the specific properties of
the topological class.
Different phenomena can occur.
Here are the known phenomena which led us to state the conjecture:
\begin{enumerate}
 \item The same asymptotic convergence for the whole topological class.
     In this case the topological class coincides with the rigidity class, see
     Example~\ref{standard}.
 \item Finitely many geometrical invariants which describe the asymptotics of
     the renormalization.
     In this case the rigidity classes foliate the topological class, see
     Examples \ref{flat1} and~\ref{cob}.
 \item Existence of open parts of the topological class in which the
     asymptotics of renormalization is determined by finitely many
     geometrical invariants.
     This happens in the coexistence phenomenon in Example~\ref{lorenz1}.
     The number and the type of the geometrical invariants is not necessarily
     the same in different open parts.
     In this case the topological class is stratified by rigidity classes.
     The boundaries of the open parts also form strata.
     See also Example~\ref{lorenz2}.
 \item Convergence of the renormalization only occurs at certain points of
     the attractor.
     This is a purely higher-dimensional phenomenon which
     gives rise to probabilistic rigidity, see Example~\ref{henon1}.
\end{enumerate}

\section{Examples}\label{examples}

Here we list the examples that were used to formulate the Rigidity Conjecture.

\subsection{Classical cases} \label{standard}

Rigidity is known to hold for many classes of systems with bounded combinatorics
and sufficient smoothness.
For example: Kleinian groups \cite{Mostow}, circle diffeomorphisms
\cites{H79,Y84}, critical circle homeomorphisms
\cites{dFdM1, Yam, Yam1, GdM, GMdM15}, unimodal maps
\cites{Lan, Sul, McM, Lyu1, Lyu2, dMP99, dFdMP},
circle maps with breakpoints \cites{KT1, KK14}.

\subsection{Flat circle maps I} \label{flat1}

For circle maps with a flat interval, critical exponent $l < 2$ and rotation
number of Fibonacci type the topological class is a codimension-$1$ manifold
and it is foliated by rigidity classes
which are codimension-$3$ submanifolds \cite{LM}.
Three geometrical invariants describe the leaves.
Invariant Cantor sets of unimodal maps of Fibonacci type with critical
exponent $l=2$ have one characterizing geometrical invariant
\cite{MilnorLyubic}.

\subsection{Flat circle maps II} \label{flat2}

Let $f$ and $g$ be two circle maps with a flat interval having critical
exponents $l_1 > l_2 > 2$ and the same rotation number of bounded type.
Both maps have a priori bounds \cite{5aut}.
A priori bounds imply that $f$ and $g$ are quasi-symmetrically conjugate on
their attractors \cite{liv}.
However, $f$ and $g$ are not smoothly conjugate since their critical exponents
differ.

\subsection{Lorenz maps I} \label{lorenz1}

Fix a critical exponent $l > 1$.
Lorenz maps of sufficiently high monotone combinatorial type exhibit the
coexistence phenomenon \cite{MW16}.
That is, there are maps $f$ and $g$ of the same topological type such that $f$
has bounded geometry and $g$ has degenerate geometry.
In particular, $f$ and $g$ are not quasi-symmetrically conjugate despite having
the same topological type and the same critical exponent.
It is conjectured that there is a codimension-$1$ stable manifold inside the
topological class and that the rest of the topological class is foliated by
finite codimension rigidity classes.
The codimensions may vary.

\subsection{Lorenz maps II} \label{lorenz2}

Fix a critical exponent $l > 1$.
Based on numerical experiments it is conjectured that the Lorenz operator has
both a fixed point and a (strict) period-two point for some low monotone
combinatorics \cite{B}.
The fixed point and the period-two point have a priori bounds and they are
quasi-symmetrically conjugate, but they are not smoothly conjugate.
The complement of the stable sets of the fixed point and the period-two point
is laminated by finite codimension rigidity classes.

\subsection{H\'enon maps I} \label{henon}

Infinitely renormalizable period-doubling H\'enon maps with different average
Jacobian are not smoothly conjugated on their attractors.
These maps are rigid when the average Jacobian $b=0$ (the unimodal case) and
when $b=1$ (the area-preserving case) \cites{CLM05, GJM}.
It is conjectured that period-doubling H\'enon maps are foliated by
codimension-$1$ rigidity classes determined by the average Jacobian $b$.
It is know that the conjugacy classes of infinitely renormalizable
period-doubling H\'enon maps are not of finite codimension.
In particular there are no finite-dimensional families which intersect all
conjugacy classes of infinitely renormalizable period-doubling H\'enon maps
\cite{HMT}.
The essential part of the dynamics of these maps is the Cantor attractor which
exists for infinitely renormalizable maps.
Hence the topologically relevant property is infinite renormalizability.
The conjugacy classes do not play a topological role in this setting.
The infinitely renormalizable maps form a codimension-$1$ manifold while the
conjugacy classes are very small.

\subsection{H\'enon maps II} \label{henon1}

The Cantor attractor of an infinitely renormalizable period-doubling H\'enon
map is uniquely ergodic.
The conjugacy between the Cantor attractors of two such maps is almost
everywhere $C^{1+\alpha}$ which respect to the unique invariant
measure, for some $\alpha>0$ \cite{ML}.
Recall that if the average Jacobians are different then the maps are not
smoothly conjugate.
This phenomenon is called probabilistic rigidity.

\subsection{Affine interval exchange transformations} \label{cob}

Almost all topological classes (which are the conjugacy classes) in the space
of affine interval exchange transformations is foliated by rigidity classes
which are one-dimensional subspaces \cite{cobo}.

\section{Observations on quasi-symmetry and a priori bounds}

Let us close with the following observations.
A system has a priori bounds if the successive renormalizations do not degenerate.
In the context of bounded combinatorics and a priori bounds in one-dimensional dynamics, the topological
class coincides with the quasi-symmetric rigidity classes.
The conjugacy between two maps is quasi-symmetric on the attractors.
This played a crucial role in the study of the convergence of renormalization.
However, quasi-symmetry does not characterize geometry.
The simplest indication of this is when two equivalent systems have both a priori bounds but different critical exponent. This is illustrated in Example~\ref{flat2} where changing the critical
exponent does not break quasi-symmetry even though the geometry changes.
A more subtle reason is given by Example~\ref{lorenz2} where both the
period-two point and the fixed point of renormalization have a priori bounds.
Hence there is quasi-symmetry but the geometry of these points cannot coincide
because of their different renormalization periods.

\bigskip

{\bf Acknowledgements:} The authors would like to thank the referee for the careful reading and the valuable remarks which helped to improve the manuscript.

\end{document}